\documentclass[10pt,onecolumn]{article}
\newcommand {\be}{\begin{equation}}
\newcommand {\ee}{\end{equation}}
\newcommand {\bey}{\begin{eqnarray}}
\newcommand {\eey}{\end{eqnarray}}
\usepackage{epsfig}
\usepackage{amsfonts}
\usepackage{amsmath}
\usepackage{mathbbol} 
\usepackage[plain]{algorithm}
\usepackage{algorithmic}

\newtheorem{theorem}{Theorem}[section]
\newtheorem{lemma}[theorem]{Lemma}

\newcommand{\qed}{\nobreak \ifvmode \relax \else
      \ifdim\lastskip<1.5em \hskip-\lastskip
      \hskip1.5em plus0em minus0.5em \fi \nobreak
      \vrule height0.75em width0.5em depth0.25em\fi}
\begin{document}

\title{Output-sensitive algorithm for generating the flats of a matroid}

\author{Alberto Montina  \\
Perimeter Institute for Theoretical Physics, 31 Caroline 
Street North, \\ Waterloo, Ontario N2L 2Y5, Canada}

\date{\today}
\maketitle

\begin{abstract}
We present an output-sensitive algorithm for generating the whole set of flats of a finite 
matroid. Given a procedure, $P$, that decides in $S_P$ time steps if a set is independent, 
the time complexity of the algorithm is $O(N^2 M S_P)$, where $N$ and $M$ are the input and 
output size, respectively. In the case of vectorial matroids, a specific algorithm is reported 
whose time complexity is equal to $O(N^2 M d^2)$, $d$ being the rank of the matroid.
In some cases this algorithm can provide an efficient method for computing zonotopes in 
$H$-representation, given their representation in terms of Minkowski 
sum of known segments.
\end{abstract}

\section{Introduction}

A matroid is a structure, introduced by Whitney~\cite{whitney}, providing an 
abstraction of the concept of independence that is common in different theories, 
such as linear algebra and graph theory. A finite matroid is defined as a pair
$({\bf W},{\cal I})$, where $\bf W$ is a finite set, called {\it ground set},
and ${\cal I}$ is a collection of subsets of $\bf W$, called independent sets, 
satisfying the following three properties~\cite{oxley}:

\begin{enumerate}
\item The set $\cal I$ is not empty;
\item If ${\bf A}\in\cal I$ and ${\bf B}\subseteq{\bf A}$, then ${\bf B}\in{\cal I}$
({\it hereditarity});
\item 
If ${\bf A},{\bf B}\in{\cal I}$ and $|\bf A|<|\bf B|$, then there exists an 
element $b\in \bf B\backslash A$ such that $\{b\}\cup\bf A\in{\cal I}$
({\it augmentation property}),
\end{enumerate}
$|\bf C|$ being the cardinality of set $\bf C$. With a slight abuse of notation,
hereafter we will denote by $b\cup\bf A$ the union of two sets $\bf A$ and 
$\{b\}$, the latter containing the single element $b$. In this paper we will
refer often to the concrete example of vectorial matroids, where
the elements of $\bf W$ are vectors of a vector space and 
independent sets are the linearly independent subsets of $\bf W$.

We are interested in an efficient algorithm for computing flats of a matroid.
Flats are subset of $\bf W$ whose properties provide an alternative 
axiomatization of matroids. Ordered by inclusion, they form a geometric 
lattice~\cite{oxley,lattice}. In the case of vectorial matroid, each flat 
with maximal rank and properly contained in $\bf W$ (called {\it hyperplane}) 
can be associated with two facets of a zonotope up to translations.
Zonotopes are polytopes, equivalently 
defined as Minkowski sum of segments or affine projections of cubes. 
They play an important role in several mathematical areas, such as 
hyperplane arrangements, box splines and partition functions~\cite{deconcini}. 
They could turn to be useful also in some problems of quantum information.

In general the
number of flats grows exponentially in $|{\bf W}|$. Thus, their
computation has an exponential time complexity. However, in
many practical problems the matroid has special properties
that considerably reduce the output size with respect to the general
case. For example, this occurs if the cardinality of some dependent 
subset of $\bf W$ is smaller than or equal to the matroid rank.
An algorithm whose running time depends only on
the input size does not take advantage of these special structures
and its complexity is exponential in any case. Conversely, an
output-sensitive algorithm, whose running time depends on both the 
input and output size, may require much less resources in cases of 
reduced output size. We will consider the extended notion of
polynomial complexity that accounts for this output sensitivity.
An algorithm is polynomial if its running time is polynomial
in both the input and output size. 

In this paper we will present an output-sensitive algorithm for 
computing all the flats of a matroid. Assuming that there is
a procedure that decides in $S_P$ steps if a set is independent,
the complexity of evaluating the flats is $O(N^2 M S_P)$, where $M$ 
is the number of flats and $N$ is the cardinality of $\bf W$.
The linearity in the number of flats is the significant feature
that makes the algorithm output-sensitive. In the concrete case of 
vectorial matroid, the procedure $P$ can be given for example by an 
algorithm that evaluates the rank of matrices. The overall complexity
of evaluating the flats is $O(N^2 M d^3)$. We will also provide a
specific optimization that requires an insignificant
increase of computation space and reduces the time complexity
to $O(N^2 M d^2)$.
Since the facets of a zonotope are identified by flats of a vectorial 
matroid up to translations, in some cases
our algorithm can provide an efficient method for evaluating
the $H$-representation of this particular polytope.
The paper is organized as follows. In Sec.~\ref{sec_def} we
give some definitions, like {\it rank}, {\it basis}, {\it closure} 
and {\it flat}. We show that each flat can be represented through
its bases. In sec.~\ref{sec_repres}, we introduce a total order
in the power set of the ground set $\bf W$ and use it to select
a representative basis for each flat. Some properties of this basis 
are proved and then used in Sec.~\ref{sec_algo}, where we present
the algorithm for generating the whole set of flats. 
In Sec.~\ref{sec_zono}, we show that our algorithm can provide in 
some cases an efficient method for evaluating the $H$-representation 
of a zonotope, given its representation in terms of Minkowski
sum of known segments.

\section{Definitions and basic properties}
\label{sec_def}

In order to define flats, it is useful to introduce some concepts, such
as basis, rank and closure. A {\it basis} is a maximal independent
set, that is, an independent set that is not properly contained in an 
independent set. The last axiom of matroid theory implies that all
the bases have the same cardinality, which is called {\it rank} of the matroid. 
Given a matroid $({\bf W},{\cal I})$ and any subset $\bf U$ of $\bf W$,
let $\cal J$ be the collection of subsets of $\bf U$ that are in $\cal I$,
then the pair $({\bf U},{\cal J})$ is a matroid. It is
called the {\it restriction} of $({\bf W},{\cal I})$ to $\bf U$~\cite{oxley}. 
Thus, each ${\bf U}\subseteq{\bf W}$ has a rank, denoted by $r({\bf U})$,
and a set of bases.

The {\it closure} of a set ${\bf U}\subseteq\bf W$, 
indicated with $\text{cl}({\bf U})$, contains the elements $w\in\bf W$ such
that $\bf U$ and ${\bf U}\cup w$ have the same rank, that is,
\be
\text{cl}({\bf U}):=\{w\in{\bf W}|r({\bf U}\cup w)=r({\bf U})\}.
\ee
The closure operator is idempotent, that is,
\be\label{idempot}
\text{cl}(\text{cl}({\bf U}))=\text{cl}({\bf U}),\;\;
\forall{\bf U}\subseteq{\bf W}.
\ee

A {\it flat} is a set ${\bf F}\subseteq{\bf W}$ that is equal to its closure,
that is, 
\be
{\bf F} \text{ is a flat }\Leftrightarrow {\bf F}=\text{cl}({\bf F}).
\ee
A flat of rank $n$ is called $n$-flat. The $(d-1)$-flat of a matroid
of rank $d$ is called {\it hyperplane}.
Flats are analogous to vector subspaces. Indeed, in the
case of vectorial matroids, it is possible to identify the flats 
with the subspaces linearly spanned by the vectors in the flats. 
In particular, $n$-flats correspond to $n$-dimensional subspaces.

By the idempotence property~(\ref{idempot}) we have that
\be
\forall{\bf X}\subseteq{\bf W},\;\; \text{cl}({\bf X})\text{ is a flat}.
\ee
Trivially, the closure operator, with the sets of flats as codomain,
is a surjective function, thus it is possible to represent every flat $\bf F$ 
through a subset whose closure is $\bf F$. The minimal sets
representing a flat $\bf F$ are independent sets with cardinality
equal to $r({\bf F})$, as implied by the followings.

\begin{lemma}
\label{cl_subset}
If ${\bf X}\subseteq{\bf Y}$ and $\mathrm{r}({\bf X})=
\mathrm{r}({\bf Y})$, then $\mathrm{cl}({\bf X})=\mathrm{cl}({\bf Y})$.
\end{lemma}

\begin{lemma}
\label{rank_cl}
$\forall{\bf X}\;\;
\mathrm{r}(\mathrm{cl}({\bf X}))=\mathrm{r}({\bf X})$.
\end{lemma}
As a direct consequence, we have the following lemma.

\begin{lemma}
The sets with minimal cardinality that represent a flat $\bf F$ are all 
the bases of $\bf F$.
\end{lemma}
{\bf Proof.} 
Let $\bf S$ be a basis of $\bf F$, thus, by definition of basis
and rank, ${\bf S}\subseteq{\bf F}$
and $\mathrm{r}({\bf S})=\mathrm{r}({\bf F})$. By
Lemma~\ref{cl_subset} and definition of flat we have that 
$\mathrm{cl}({\bf S})=\mathrm{cl}({\bf F})={\bf F}$, that is,
$\bf S$ represents the flat ${\bf F}$. This set is also
minimal because of Lemma~\ref{rank_cl}. Indeed, if $\bf R$
represents the flat $\bf F$, then by definition
$\mathrm{cl}({\bf R})={\bf F}$. By Lemma~\ref{rank_cl},
this implies that
$\mathrm{r}({\bf R})=\mathrm{r}({\bf F})$,
that is, the cardinality of $\bf R$ is not smaller than
the cardinality of a basis of $\bf F$.

Thus, flats can also be identified with the class of their bases, just 
as in linear algebra a set of $k$ independent vectors identifies
a $k$-dimensional linear subspace. This provides a simplification 
in the representation of flats, since it is not necessary to enumerate 
the whole set of its elements. The mapping from the bases to the flats 
is surjective, but in general is not bijective.
In the next section, we introduce a rule for
selecting a representative basis that will turn to be fundamental
for developing our algorithm.

\section{Denoting flats through representative bases}
\label{sec_repres}
To select a representative basis of a flat, we introduce a total 
order on the subsets of $\bf W$ and associate each flat with its 
{\it first} basis (``first'' with respect to the total order).
The order is defined as follows.

Given a matroid $({\bf W},{\cal I})$, we order
the elements of ${\bf W}$ by appending an integer
$i\in[1,N]$ to each element $w_i\in\bf W$, $N$ being
the cardinality of $\bf W$. Then we represent a subset 
${\bf U}\subseteq{\bf W}$ through a $N$-digit binary number
by setting the $i$-th digit equal to $1(0)$ if $w_i$
is (not) an element of $\bf U$, for every $i\in[1,N]$.
In other words, given any collection
$\{b_1,...,b_k\}$ of indices with $b_{n-1}<b_n$, we label 
the subset
${\bf X}=\{w_{b_1}, w_{b_2},...,w_{b_k}\}
\subseteq\bf W$ with the number 
\be
L({\bf X})=\sum_{n=1}^k 2^{b_n-1}.
\ee
For example, 
$$
L(\{w_1,w_3,w_7\})=1000101_2,
$$
which is equal to $69_{10}$ in decimal basis. 
We call the most significant nonzero bit of a binary number $b$
{\it leading} digit of $b$. 

By attaching the label $L$ to each subset of $\bf W$, we have 
introduced a total order in the power set of $\bf W$. We associate 
each flat $\bf F$ with the basis that has the smallest label $L$,
which we call {\it pointer} of $\bf F$. In particular an $n$-pointer 
is the pointer of an $n$-flat. The pointer of $\bf F$
will be indicated with $p({\bf F})$. It is important to
distinguish the pointer of a flat from its label $L$.
The pointer of $\bf F$ is the label of its first basis,
which in general is different from the label of $\bf F$,
unless $\bf F$ is independent. In binary representation, the number 
of non-zero digits of the pointer and label of a flat $\bf F$ is 
equal to $\mathrm{r}({\bf F})$ and $|{\bf F}|$, respectively.

\begin{theorem}
\label{del_bit}
Let $p_0$ be the pointer of a $n$-flat in binary representation, then 
the number $s$ obtained from $p_0$ by replacing the leading digit
with $0$ is the pointer of an $(n-1)$-flat.
\end{theorem}

In order to prove it, we need the following property.
\begin{lemma}
\label{add_vect}
Let ${\bf I}_1$ and ${\bf I}_2$ be two independent subset of 
$\bf W$ such that $\mathrm{cl}({\bf I}_1)=\mathrm{cl}({\bf I}_2)$,
then, for every $w\in{\bf W}$, 
$\mathrm{cl}({\bf I}_1\cup w)=\mathrm{cl}({\bf I}_2\cup w)$.
\end{lemma}
{\bf Proof.}
If $w\in\mathrm{cl}({\bf I}_1)=\mathrm{cl}({\bf I}_2)$, then
the lemma is a direct consequence of the implication 
$w\in\mathrm{cl}({\bf X})\Rightarrow\mathrm{cl}({\bf X})=\mathrm{cl}
({\bf X}\cup w)$, which can be easily obtained from 
lemma~\ref{cl_subset}. Thus, let us consider the case 
$w\notin\mathrm{cl}({\bf I}_1)=\mathrm{cl}({\bf I}_2)$,
that is, we assume that ${\bf I}_1\cup w$ and ${\bf I}_2\cup w$
are independent. We have to prove that if an element
$v\notin\mathrm{cl}({\bf I}_1\cup w)$, then
$v\notin\mathrm{cl}({\bf I}_2\cup w)$ and vice versa.
Suppose that $v\notin\mathrm{cl}({\bf I}_1\cup w)$,
then ${\bf I}_1\cup w\cup v$ is independent. Since also
${\bf I}_2\cup w$ is independent and 
$|{\bf I}_1\cup w\cup v|>|{\bf I}_2\cup w|$, 
by axiom 3 of matroid theory (augmentation property) there is an element
$b\in{\bf I}_1\cup w\cup v$ such that 
${\bf I}_2\cup  w\cup  b$ is independent. Clearly, $  b$
is not in ${\bf I}_1$, since $\mathrm{cl}({\bf I}_1)=\mathrm{cl}({\bf I}_2)$
and $  b\notin{\bf I}_2$.
Furthermore, $  b$ cannot be equal to $  w$, thus $  b=  v$.
This implies that ${\bf I}_2\cup  w\cup  v$ is independent,
that is, $  v\notin\mathrm{cl}({\bf I}_2\cup  w)$.
Also the inverse implication is true.
Thus, every element that is not in $\mathrm{cl}({\bf I}_1\cup  w)$
is not in $\mathrm{cl}({\bf I}_2\cup  w)$ and
vice versa, that is, $\mathrm{cl}({\bf I}_1\cup  w)$ and
$\mathrm{cl}({\bf I}_2\cup  w)$ are equal. $\square$
\newline

{\bf Proof of Theorem~\ref{del_bit}.} 
Let ${\bf I}_1$ be the independent set with label $s$.
Suppose that $s$ is not a pointer, thus there is an independent set
${\bf I}_2$ with $L({\bf I}_2)<L({\bf I}_1)$  such that 
$\mathrm{cl}({\bf I}_2)=\mathrm{cl}({\bf I}_1)$.
Let $k$ be the position of the leading digit of $p_0$, 
then, by definition of $s$ and ${\bf I}_1$,
${\bf I}_1\cup  w_k$ is the independent set
pointed to by $p_0$. By lemma~\ref{add_vect}, both  
${\bf I}_1\cup  w_k$ and  ${\bf I}_2\cup  w_k$
are bases of the same flat and furthermore
$L({\bf I}_2\cup  w_k)<L({\bf I}_1\cup  w_k)$, since 
$L({\bf I}_2)<L({\bf I}_1)$, but this is impossible because
$p_0=L({\bf I}_1\cup  w_k)$ is a pointer.
$\square$ \newline

This lemma implies that each $(n+1)$-pointer can be generated from 
some $n$-pointer by setting one of the digits at the left of the
leading digit equal to $1$. For example, if $a=10010011_2$ is a $4$-pointer, 
then the number $b=00010011_2$ is a $3$-pointer for Theorem~\ref{del_bit}.
The $4$-pointer $a$ is generated from the $3$-pointer $b$ by replacing
the $8$-th zero digit of $b$ with $1$. Thus, given a collection of 
$n$-pointers, this procedure of replacement generates a set of 
labels that contains the set of all the 
$(n+1)$-pointers. In general the inclusion is strict, that is, the 
procedure of adding a bit $1$ to an $n$-pointer does not necessarily 
give an $(n+1)$-pointer and we need a criterion for discarding
labels that are not pointers.

\begin{theorem}
\label{cond_pointer}
Let `$s$' and `$\delta$' be the label of an subset ${\bf X}$ and the position of 
the leading digit of $s$, respectively. Let ${\bf Y}_i$ be the set obtained 
from $\bf X$ by removing the element $  w_j\in{\bf X}$ with $j>i$. The 
integer $s$ is a pointer if and only if $\bf X$ is independent and,
for every $  w_k\notin{\bf X}$ with $k<\delta$, 
$  w_k\notin\mathrm{cl}({\bf X})$ or $  w_k\in\mathrm{cl}({\bf Y}_k)$.
\end{theorem}

Note that ${\bf Y}_i$ and ${\bf Y}_j$ are not necessarily different
if $i\ne j$.
In order to prove this lemma, we need three properties. The first one is 
known as the {\it Mac Lane-Steinitz exchange property}~\cite{oxley}.

\begin{lemma}
\label{maclane}
Given a subset ${\bf X}\subseteq{\bf W}$ and an element $  w\in{\bf X}$, 
if $  v\in\mathrm{cl}({\bf X})$ and 
$  v\notin\mathrm{cl}({\bf X}\backslash  w)$, then 
$\mathrm{cl}({\bf X}\backslash  w\cup  v)=\mathrm{cl}({\bf X})$.
\end{lemma}
The second one is the {\it basis exchange property}~\cite{oxley}.
\begin{lemma}
\label{basis_exchange}
If ${\bf B}_1$ and ${\bf B}_2$ are two bases of a subset of $\bf W$
and $  w\in{\bf B}_1$, then there is an element $  v\in
{\bf B}_2\backslash{\bf B}_1$ such that 
${\bf B}_1\backslash  w\cup  v$ is a basis.
\end{lemma}
Finally, the last property is a consequence of the augmentation axiom.
\begin{lemma}
\label{deletion}
Let ${\bf Y}$ be a subset of ${\bf X}\in{\cal I}$. If 
$  w\notin\mathrm{cl}({\bf Y})$, then there is an
element $  v\in{\bf X}\backslash{\bf Y}$ such that $  w\notin
\mathrm{cl}({\bf X}\backslash  v)$.
\end{lemma}
Before proving theorem~\ref{cond_pointer}, let us prove the last
property.

{\bf Proof of lemma~\ref{deletion}.} First, we assume that
$  w\in{\bf X}$. It is is clear that if $  v=  w$
then $  w\notin\mathrm{cl}({\bf X}\backslash  v)$, since
$\bf X$ is independent. Furthermore $  v\in{\bf X}\backslash{\bf Y}$,
since $  w\in{\bf X}$ and $  w\notin\backslash{\bf Y}$,
and the conclusion of the lemma is proved.
Now we assume that $  w\notin{\bf X}$. 
For the augmentation property, it is possible to construct an
independent set by adding
$|{\bf X}|-|{\bf Y}\cup  w|$ elements in ${\bf X}\backslash{\bf Y}$ 
to ${\bf Y}\cup  w$. The obtained set is equal to
${\bf X}\cup  w$ minus some element $  v\in{\bf X}\backslash{\bf Y}$. 
Thus, ${\bf X}\backslash  v\cup  w$ is independent, that is,
$  w\notin\mathrm{cl}({\bf X}\backslash  v)$. $\square$

{\bf Proof of theorem~\ref{cond_pointer}.} First we prove one direction of 
the implication and assume that $s$ is the pointer of a flat $\bf F$.
By definition $\bf X$ is independent. Suppose that the other part of
the conclusion is false,
thus there is an element $  w_k\notin{\bf X}$ with $k<\delta$ such that 
$  w_k\in\mathrm{cl}({\bf X})$ and $  w_k\notin\mathrm{cl}({\bf Y}_k)$. 
Last condition and lemma~\ref{deletion} imply that 
there is an element $  w_l\in{\bf X}\backslash{\bf Y}_k$ such that
$  w_k\notin\mathrm{cl}({\bf X}\backslash  w_l)$. By definition
of ${\bf Y}_k$ we have that
$l>k$. Thus, since $  w_k\in\mathrm{cl}({\bf X})$ and
$  w_k\notin\mathrm{cl}({\bf X}\backslash  w_l)$, 
by lemma~\ref{maclane} we have that
$\mathrm{cl}({\bf X}\backslash  w_l\cup  w_k)=\mathrm{cl}({\bf X})$,
that is, ${\bf X}\backslash  w_l\cup  w_k$ is a basis of
$\bf F$, but this is impossible because 
$L({\bf X}\backslash  w_l\cup  w_k)<s$ ($l$ is greater than $k$) 
and $s$ is a pointer,
thus one direction of the implication is proved. 

Let us prove the other direction. Suppose that $s$ is not a pointer,
then there is a number $s_1<s$ such that $L^{-1}(s_1)\equiv\bar{\bf X}$ 
and $L^{-1}(s)={\bf X}$ are bases of the same flat $\bf F$. Let $  w_l$ 
be the element with largest subscript $l$ such that $  w_l\in{\bf X}$ and 
$  w_l\notin\bar{\bf X}$. Since $s_1<s$, this element exists. 
By lemma~\ref{basis_exchange} 
there is an element $  w_k\in\bar{\bf X}\backslash\bf X$,
such that ${\bf X}\backslash  w_l\cup  w_k$ is independent, that
is, $  w_k\notin\mathrm{cl}({\bf X}\backslash  w_l)$.
The elements with subscript larger than $l$ are in
$\bf X$ if and only if are in $\bar{\bf X}$. Since $  w_k\notin{\bf X}$
and $  w_k\in\bar{\bf X}$, then $k$ cannot be larger than $l$, thus
we have that $k<l$. Because of this inequality, ${\bf Y}_k$ is a subset 
of ${\bf X}\backslash  w_l$. The relations 
${\bf Y}_k\subseteq{\bf X}\backslash  w_l$ and
$  w_k\notin\mathrm{cl}({\bf X}\backslash w_l)$ imply that
$  w_k\notin\mathrm{cl}({\bf Y}_k)$. It is also clear that 
$  w_k\in\mathrm{cl}({\bf X})$, since $  w_k$
is in $\bar{\bf X}$ and $\mathrm{cl}(\bar{\bf X})=\mathrm{cl}({\bf X})$.   
$\square$

\section{Algorithm for generating the flats}
\label{sec_algo}

Theorems~\ref{del_bit} and \ref{cond_pointer} are the two key ingredients
of our algorithm for calculating the flats of a matroid. The idea
is generating recursively the $i$-pointers from the lower-dimensional
$(i-1)$-pointers. The flats are then generated from their pointers.
More precisely, a set of labels is generated from each $(i-1)$-pointer
by setting one of the digits at the left of the leading digit equal to $1$.
If $l$ is the position of the leading digit of an $(i-1)$-pointer and $N$ 
is the number of elements in $\bf W$, then a $(i-1)$-pointer generates 
$N-l$ labels. The set of labels generated from all the $(i-1)$-pointers
contains the whole set of $i$-pointers. Theorem~\ref{cond_pointer} provides
an efficient method for discarding labels that are not pointers.

Since the structure of flats is unaffected by
the presence of loops and parallel elements, we will assume
without loss of generality that they are absent, that it,
we assume that the matroid is simple~\cite{oxley}. Loops
are dependent subsets of $\bf W$ with cardinality equal to $1$.
In the case of vectorial matroids, a loop is a zero vector.
Parallel elements are pairwise dependent vectors.
Denoting by $N$ and $d$ the cardinality and the rank of the matroid, 
respectively, the algorithm for generating the pointers is as follows.

\begin{algorithm}[H]
%\caption{Generation of the pointers}
\parbox[t]{0.5in}{\bf Input:} The set of $N$ $1$-pointers \;\{the elements in ${\bf W}$\}
\begin{algorithmic}[1]
\STATE Set $\bar M$ equal to number of $1$-pointers. \COMMENT{$:=N$}
\FOR{$i=2$,...,$d-1$} 
\FOR{$j=1$,...,$\bar M$}
\STATE Set $l$ equal to the position of the leading digit of the $j$-th $(i-1)$-pointer.
\FOR{$\delta=l+1$,...,$N$}
\label{loop_l}
\STATE Generate from $j$-th $(i-1)$-pointer a label $s$ by setting the $\delta$-th digit equal 
to $1$.
\IF{$s$ is a pointer  \COMMENT{This is checked through theorem~\ref{cond_pointer}}}
\label{check_line} 
\STATE Store $s$ as a new $i$-pointer.
\ENDIF
\ENDFOR
\ENDFOR
\STATE Set $\bar M$ equal to the number of $i$-pointers.
\ENDFOR
\end{algorithmic}
\parbox[t]{0.6in}{\bf Output:} The whole set of pointers
\end{algorithm}

The check at line~\ref{check_line} is performed through theorem~\ref{cond_pointer}. Thus, 
it requires to verify that the set $\bf X$ pointed to by $s$ is independent and, for 
every $  w_k\notin{\bf X}$ with $k<\delta$, the set ${\bf X}\cup  w_k$ is independent 
or ${\bf Y}_k\cup  w_k$ is dependent. Given a procedure, $P$, that decides
with $S_P$ steps if a set is independent, the time complexity of the check at 
line \ref{check_line} is $\mathrm{O}(N S_P)$. The overall time complexity of the algorithm is
$O(N^2 M S_P)$, where $M$ is the total number of flats, that is, $M=\sum_{i=1}^{d-1}M_i$,
$M_i$ being the number of $i$-flats.

In the specific case of a vectorial matroid, the procedure $P$ can be provided
for example by a routine that evaluates the rank of matrices. Indeed any set 
${\bf X}\subseteq{\bf W}$ 
of a vectorial matroid can be seen as a $d\times|{\bf X}|$ matrix (called column 
matroid), the columns being the elements of $\bf X$. The rank of $\bf X$ is 
the rank of the matrix. The set is independent if the rank is equal
to $|{\bf X}|$. The time complexity of evaluating the flats by
using this routine is $O(N^2 M d^3)$. In this scheme the rank of the matrices
is evaluated without taking advantage of the similarity of their structure.
Indeed it is possible to reduce the time complexity
by a slight increase of the computational space that exploits this similarity.
Let $p_0$ and ${\bf Z}$ be an $(i-1)$-pointer and its correponding independent
set. $l$ is the position of the leading digit of $p_0$. $N-l$ labels 
are generated from $p_0$ by setting the $\delta$-th digit equal to $1$, where $\delta$ is an 
integer ranging between $l+1$ and $N$ (see line~\ref{loop_l} in the algorithm).
We denote by ${\bf X}^{(\delta)}$ the sets associated with 
the generated labels. In line~\ref{check_line} of the algorithm, first we have
to verify that ${\bf X}^{(\delta)}$ is independent for each $\delta$. Since sets
with different values of $\delta$ differ in one element and share the subset $\bf Z$,
the best strategy to decide if the sets ${\bf X}^{(\delta)}$ are independent is, first,
reducing $\bf Z$ to the row echelon form using row operations and,then, performing
the same raw operations on the added column in each ${\bf X}^{(\delta)}$.
The corresponding time complexity is $O(d^2 N)$, taking into account that
$d\le N$. A similar strategy can be used also in checking the independence
of ${\bf X}\cup  w_k$ and the dependence of ${\bf Y}_k\cup  w_k$.
In this way it is possible to reduce the overall complexity of generating the flats
to $O(N^2 M d^2)$.

\section{Minkowski sum of segments: zonotope}
\label{sec_zono}

The algorithm for the computation of flats can be useful in some cases
for calculating the $H$-representation of a zonotope when it is 
represented as Minkowski sum of known segments.
The Minkowski sum of two sets $A$ and $B$ in a vector space is the
set obtained by adding every vector of $A$ to every vector of $B$,
that is,
\be
A+B=\{\vec a+\vec b|\vec a\in A,\vec b\in B\}.
\ee
The zonotope is a polytope defined as the Minkowski sum of segments.
Up to a translation, it is the set of vectors
\be
\vec v=\sum_{k=0}^M \lambda_k \vec w_k,
\ee
where $\lambda_k\in[0:1]$. Each vector $\vec w_k$ and the zero vector 
$\vec 0$ are the two vertices of each segment summed up. 
The set ${\bf W}=\{\vec w_k|k\in[1:M]\}$ is the ground set of a vectorial
matroid, whose independent sets are the sets of linearly independent
vectors.
Let $d$ be the dimension of the zonotope, that is, the maximal
number of independent vectors in $\bf W$.
Without loss of generality, we assume that $d$ is also the
dimension of the vector space.

Each facet is parallel to a $(d-1)$-flat of the matroid, thus
its normal vector is orthogonal to any basis of the $(d-1)$-flat. 
The computation of each normal vector by a set of $d-1$ vectors
has a complexity that scales like $d^3$. We denote a normal
vector with $\vec n_i$, where the subscript $i$ is an integer
that goes from $1$ to $M_{d-1}$, $M_{d-1}$ being the number
of $(d-1)$-flats.
It can be proved that for each vector $\vec n_i$ there exist two facets 
defined by the inequalities
\be
\vec n_i\cdot\vec x\le\sum_k \theta(\vec n_i\cdot\vec w_k)\vec n_i\cdot\vec w_k
\ee
and
\be
\vec n_i\cdot\vec x\ge-\sum_k\theta(-\vec n_i\cdot\vec w_k) \vec n_i\cdot\vec w_k,
\ee
where $\theta(x)$ is the step function $\theta(x>0)=1$, $\theta(x<0)=0$.
These inequalities define the zonotope in $H$-representation.
Thus, the computation of the zonotope in $H$-representation is achieved
by evaluating the $(d-1)$-flats of the matroid $\bf W$. 

The output in this problem is the set of half-planes, thus its
size is $M_{d-1}$. In general our algorithm for the computation of
the flats does not allows us to solve this problem in polynomial
time with respect to the output size, since the algorithm presented in
the previous section is linear in the total number of flats $M$ and
in the worse case $M$ could be exponentially greater than $M_{d-1}$. However,
in many practical problems $M$ can scale linearly in $M_{d-1}$
and the input size. Suppose for example that the vectors $\vec w_k$ are
in general position and the matroid rank $d$ is smaller than $N/2$.
The number of $k$-flats is $\frac{N!}{k!(N-k)!}$ and grows monotonically
in $k(<d-1)$. This implies that $M$ scales at most like $d M_{d-1}$.
This linear scaling can be present also in the case of special structure
for which an output-sensitive algorithm provides an advantage.
It is worthwhile to note that the best algorithm for the evaluation of 
a zonotope in $H$-representation has a complexity that is quadratic in 
the output size in any case~\cite{seymour}. In the subclass of problems 
where the number of overall flats is a linear function of the number of 
hyperplanes, our method for the generation of the zonotope in 
$H$-representation is linear in the output size. An open question is
determining how much large is this subclass

\section*{Acknowledgments}
\nonumber

The Author acknowledges useful discussions with Komei Fukuda,
Hans Raj Tiwary and Erik Schnetter.
Research at Perimeter Institute for Theoretical Physics is
supported in part by the Government of Canada through NSERC
and by the Province of Ontario through MRI.

\end{document}